\documentclass{amsart}
\usepackage{graphicx}
\usepackage{newlattice}

\newcommand{\cll}{C_\tup{ll}}
\newcommand{\cul}{C_\tup{ul}}
\newcommand{\clr}{C_\tup{lr}}
\newcommand{\cur}{C_\tup{ur}}

\theoremstyle{plain}
\newtheorem{theorem}{Theorem}

\begin{document}
\title[The congruence representation theorem for semimodular lattices]{A short proof\\ of the congruence representation theorem\\ for semimodular lattices}  
\author{G. Gr\"{a}tzer} 
\address{Department of Mathematics\\
  University of Manitoba\\
  Winnipeg, MB R3T 2N2\\
  Canada}
\email[G. Gr\"atzer]{gratzer@me.com}
\urladdr[G. Gr\"atzer]{http://server.maths.umanitoba.ca/homepages/gratzer/}

\author{E.\,T. Schmidt}
\address{Mathematical Institute of the Budapest University of
        Technology and Economics\\  
        H-1521 Budapest\\
        Hungary}
\email[E.\,T. Schmidt]{schmidt@math.bme.hu}
\urladdr[E.\,T. Schmidt]{http://www.math.bme.hu/\~{}schmidt/}

\date{March 30, 2013}

\thanks{The second author was supported by
the Hungarian National Foundation for Scientific Research (OTKA),
grant no.\ K77432.}
\begin{abstract}
In a 1998 paper with H. Lakser, 
the authors proved that every finite distributive lattice $D$
can be represented as the congruence lattice 
of a~finite \emph{semimodular lattice}.

Some ten years later, the first author and E. Knapp proved a much stronger result, proving the representation theorem for 
\emph{rectangular lattices}.

In this note we present a short proof of these results.
\end{abstract}

\maketitle

\section{Introduction}\label{S:Intro}
In \cite{GLS98a}, the authors with H. Lakser
proved the following result:

\begin{theorem}\label{T:semimod}
Let $D$ be a finite distributive lattice. 
Then there is a \emph{planar semimodular lattice}~$K$ such that
\[
   D \iso \Con K.
\]
\end{theorem}

A stronger result was proved some 10 years later. 
To state it, we need a few concepts.

Let $A$ be a planar lattice. 
A \emph{left corner} (resp., \emph{right corner}) of the lattice $A$ 
is a doubly-irreducible element in $A - \set{0,1}$ 
on the left (resp., right) boundary of~$A$. 

We define a \emph{rectangular lattice} $L$, 
as in G.~Gr\"atzer and E. Knapp \cite{GK09}, 
as a planar semimodular lattice 
that has exactly one left corner, $u_l$, 
and exactly one right corner, $u_r$, 
and they are complementary---that is, 
$u_l \jj u_r = 1$ and $u_l \mm u_r = 0$. 

The first author and E. Knapp \cite{GK09} proved 
the following much stronger form of Theorem~\ref{T:semimod}:

\begin{theorem}\label{T:rectangular}
Let $D$ be a finite distributive lattice. 
Then there is a \emph{rectangular lattice}~$K$ such that
\[
   D \iso \Con K.
\]
\end{theorem}

In this note we present a short proof of this result.

\section{Notation}\label{S:Preliminaries}
We use the standard notation, see \cite{LTF}.

For a rectangular lattice $L$, 
we use the notation $\cll = \id{u_l}$,  $\cul = \fil{u_l}$,  
$\clr = \id{u_r}$,  $\cur = \fil{u_r}$ 
for the four boundary chains; 
if we have to specify the lattice $L$, 
we write  $\cll(L)$, and so on.
(See G. Cz\'edli and G. Gr\"atzer \cite{CGa} 
for a survey of semimodular lattices, in general, 
and rectangular lattices, in particular.)

\begin{figure}[h!]
\centerline{\includegraphics[scale=1]{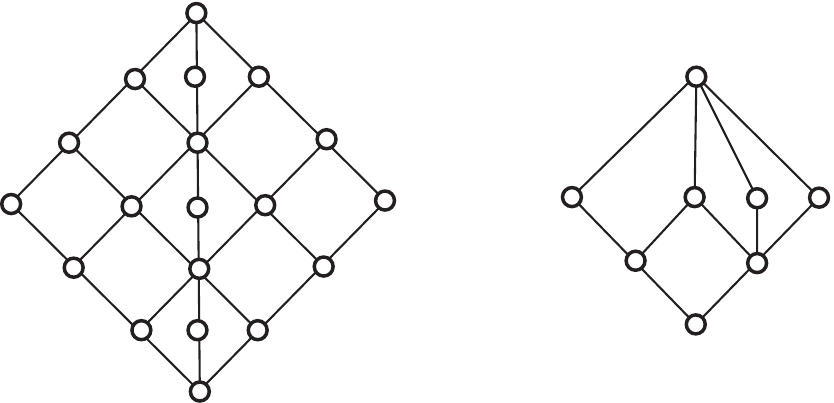}}
\caption{The $\SM 3$-grid for $n=3$ and the lattice $\SaS 8$}\label{F:grid+S8}

\medskip

\medskip

\medskip

\medskip

\medskip

\centerline{\includegraphics[scale=1.0]{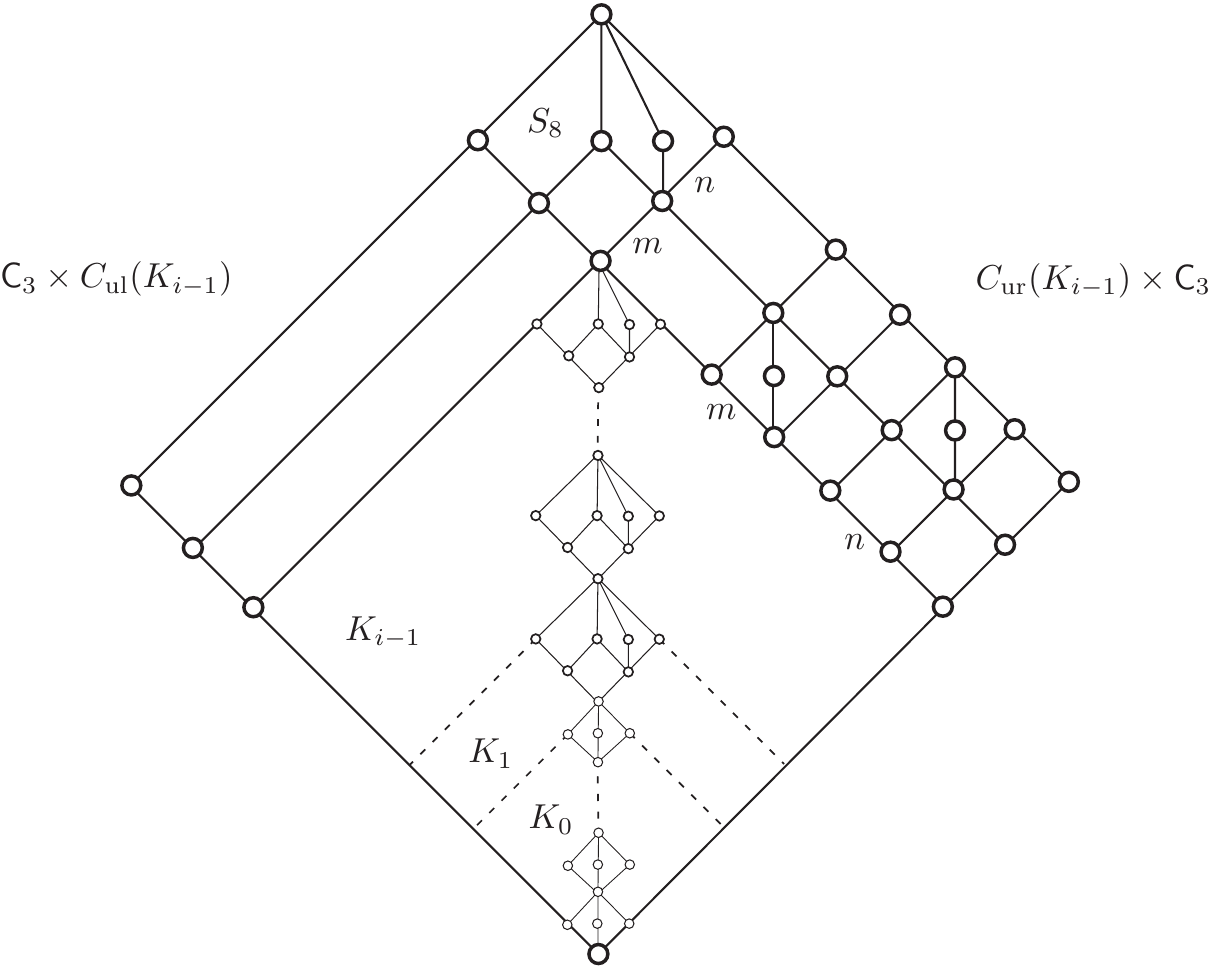}}
\caption{A sketch of the lattice $K_i$ 
for $n \geq 3$ and $3 < i \leq e$}\label{F:K}
\end{figure}

\section{Proof}\label{S:Proof}
Let $D$ be the finite distributive lattice of Theorem~\ref{T:rectangular}.
Let $P = \Ji D$. Let $n$ be the number of elements in $P$
and $e$ the number of coverings in $P$.

We shall construct a rectangular lattice $K$
representing $D$ by induction on $e$. 
Let $m_i \prec n_i$, for $1 \leq i \leq e$, list all coverings of $P$.
Let $P_j$, for $0 \leq j \leq e$, be the order we get from $P$ by removing the coverings
$m_i \prec n_i$ for $j < i \leq e$. 
Then $P_0$ is an antichain and $P_e = P$.

For all $0 \leq i \leq e$,
we construct a rectangular lattice $K_i$ inductively.
Let $K_0 = \SC {n+1}^2$ be a grid,
in which we replace the covering squares of the main diagonal
by covering $\SM 3$-s; see Figure~\ref{F:grid+S8} for $n = 3$.
Clearly, this lattice is rectangular 
and $\Con K_0$ is the boolean lattice with $n$ atoms.

Now assume that $K_{i-1}$ has been constructed. 
Let the three-element chain $0 \prec m_i \prec n_i$ 
be represented by the lattice $\SaS 8$, 
see Figure~\ref{F:grid+S8}. 

Take the four lattices 
\[
   \SaS 8,\ K_{i-1},\ \SC 3 \times \cul(K_{i-1}), 
   \cur(K_{i-1}) \times \SC 3
\]
and put them together as in Figure~\ref{F:K}, 
where we sketch $K_{i-1}$ for $n \geq 3$ and $3 < i \leq e$.
We add two more elements to turn two covering squares
into covering $\SM 3$-s, see Figure~\ref{F:K},
so that the prime interval of $\SaS 8$ 
marked by $m$ defines the same congruence
as the prime interval of $K_{i-1}$ marked by~$m$; 
and the same for $n$. 
Let $K_i$ be the lattice we obtain. 
The reader should have no trouble to directly verify that 
$K_i$ is a rectangular lattice. 
(See G.~Cz\'edli and G.~Gr\"atzer \cite{CGa} 
for general techniques that could be employed.)

The lattice $K$ for Theorem~\ref{T:rectangular} is the
lattice $K_{e}$.

See G. Gr\"atzer \cite{GGa} for a comparison 
how this short proof compares to the proofs in 
G. Gr\"atzer, H. Lakser, and E.\,T. Schmidt \cite{GLS98a} 
and in G. Gr\"atzer and E. Knapp~\cite{GK09}.

\end{document}